\documentclass[12pt,a4paper]{article}
\usepackage{amstext}
\usepackage{fancyhdr}
\usepackage{amsfonts,graphicx,bezier, amssymb}
\usepackage{caption}
\usepackage{etoolbox}
\usepackage{caption}

\usepackage{graphics,graphicx, bezier, float, color}
\usepackage{amsmath,amsthm,amssymb,amsfonts,amscd,enumerate, array,latexsym,epsfig,psfrag}
\newtheorem{defn}{Definition}[section]
\let\olddefinition\defn
\renewcommand{\defn}{\olddefinition\normalfont}
 \newtheorem{prop}{Proposition}[section]
 
  \newtheorem{thm}{Theorem}[section]
  \newtheorem{cor}[thm]{Corollary}

  \newtheorem{exam}{Example}[section]
\newenvironment{prf}{\noindent{\bf{Proof:}}~~}{\hfill\rule{1ex}{1ex}\vskip1.5ex}
\newcommand{\Z}{\mathbb Z}

\newcommand{\Ra}{\Rightarrow}

\newcommand{\beqa}{\begin{eqnarray}}
\newcommand{\enqa}{\end{eqnarray}}
\newcommand{\beq}{\begin{eqnarray*}}
\newcommand{\enq}{\end{eqnarray*}}

\begin{document}
\begin{center}
{\bf\Large The locally nilradical for modules over commutative rings }

\vspace*{0.5cm}

 \begin{center}

Annet Kyomuhangi and David Ssevviiri\footnote{Corresponding author}\\
Department of Mathematics, Makerere University\\
 P.O BOX 7062, Kampala Uganda\\
E-mail: annet.kyomuhangi@gmail.com and ssevviiri@cns.mak.ac.ug\\

 \end{center}

\end{center}

\vspace*{0.5cm}
\begin{abstract}
Let $R$ be a commutative unital ring and $a\in R.$ We introduce and study properties of a functor $a\Gamma_{a}(-),$ called the  locally nilradical  on the category of $R$-modules.
  $a\Gamma_{a}(-)$ is a generalisation of both the torsion functor (also called section functor) and  Baer's lower nilradical for modules. Several local-global properties of the functor $a\Gamma_{a}(-)$ are established. As an application, results about reduced $R$-modules are obtained and hitherto unknown ring theoretic radicals as well as structural theorems are deduced.

   \end{abstract}


\vspace*{0.4cm}

{\bf MSC 2010} Mathematics Subject Classification: 16S90, 16N80, 13D45

\section{Introduction}

Radicals are a good tool to study the structure of rings and modules over rings. There are several radicals in the literature about rings and modules which include among others; Baer's lower nilradical (also called the prime radical), K\"{o}the's upper nilradical, Andrunakievich's generalised nilradical (also called the completely prime radical), Jacobson radical and Brown-McCoy radical. In this article, we introduce and study a radical called the locally nilradical for modules over commutative rings. Radical theory also exists for abelian categories, and it is what is termed as torsion theory.

\paragraph\noindent
Throughout this paper, all rings $R$ are commutative and unital. The
 category of  all left $R$-modules is denoted by $R$-Mod. If $M \in R$-Mod, then $\sqrt{(0:M)}$ denotes the radical ideal of $(0:M),$ i.e., $\sqrt{(0:M)}= \left \{r\in R~|~r^{k}\in (0:M) ~\text{for some } k\in \mathbb{Z}^{+} \right \}.$
Let $M$ be an $R$-module and ${\bf{\mathfrak{a}}}$  an ideal of $R.$ The $\bf{\mathfrak{a}}$-torsion (also called the section) functor is defined by: $$\Gamma_{{\bf{\mathfrak{a}}}}:R\text{-Mod} \to R\text{-Mod}$$
$$M\mapsto \Gamma_{{\bf{\mathfrak{a}}}}(M), $$ where $\Gamma_{{\bf{\mathfrak{a}}}}(M)$ is the submodule of $M$ given by $\Gamma_{\bf{\mathfrak{a}}}(M):= \left \{ m\in M ~| ~{\bf{\mathfrak{a}}}^{k}m= 0 ~\text{for some } k\in \mathbb{Z}^{+} \right \}.$
On modules defined over Noetherian rings, this functor is left exact and a radical, see \cite{Rohrer}. Its right derived functor $H_{\mathfrak{a}}^i(-)$ is what is called the local cohomology functor with respect to $\mathfrak{a}.$ For more information about local cohomology, see \cite{Brodmann}. If $R$ is a ring, $a\in R$ and $\bf{\mathfrak{a}}$ an ideal of
$R$ generated by $a,$ then it is easy to see that $\Gamma_{\bf{\mathfrak{a}}}(M)= \Gamma_{a}(M),$ where  $\Gamma_{a}(M):= \left \{ m\in M ~| ~a^{k}m= 0 ~\text{for some } k\in \mathbb{Z}^{+} \right \}.$

\paragraph\noindent

By generalising the torsion functor, we define a new functor:

$$a\Gamma_{a}:R\text{-Mod} \to R\text{-Mod}$$

$$M\mapsto a\Gamma_{a}(M) ,$$ called the \emph{locally nilradical} which associates to every $R$-module $M$ a submodule $a\Gamma_{a}(M)$ for every $a\in R$, where
$a\Gamma_{a}(M):=\left\{am~|~a^{k}m=0, ~m\in M, \text{for some}~k\in \Z^{+} \right\}$, i.e.,
left multiplication by $a$ of the submodule $\Gamma_a(M)$.

\paragraph\noindent

$a\Gamma_{a}(M)$ is contained in the envelope $E_{M}(0)$ of $M$ which has been considered in the literature as a module analogue of the set of nilpotent elements of a ring. Secondly, we observe that if $M$ is the $R$-module $R,$ then $a\Gamma_{a}(-)$ associates to $R$ a nil ideal $a\Gamma_{a}(R)$ of $R.$ For if $x\in a\Gamma_{a}(R),$ then $x= ar $ and $a^{k}r =0$ for some $k\in \Z^{+}.$ It follows that $x^{k}= (ar)^{k}= a^{k}r^{k}=0.$ So, $x$ is nilpotent and $a\Gamma_{a}(R)$ is nil. We use the adjective $``\text{locally}"$ because $a\Gamma_{a}(-)$ gives the local behaviour for a given element $a\in R$ as opposed to the global picture which is for all $a\in R$ given by the nilradical $\mathcal{N}(R).$  To be precise, we show that $\bigcup\limits_{a\in R} a\Gamma_{a}(R)= \mathcal{N}(R).$

\paragraph\noindent
This  paper is devoted to studying properties of the functor $a\Gamma_{a}(-)$. We list some of them below.

\begin{itemize}

\item  $a\Gamma_{a}(-)$ is a radical on the category $R$-Mod, (Proposition \ref{nilrad}).

\item For a Noetherian local ring $R$ of characteristic $p,$ the Frobenius functor $F_{R}(-)$ is exact on $R$-Mod if and only if for any $a\in R$, the functor $a\Gamma_{a}(-)$ is trivial on $R$-Mod, (Theorem \ref{Frob}).

\item For any ring $R,$ and $a\in R,$ $a\Gamma_{a}(R)[x]= a\Gamma_{a}(R[x]),$ (Theorem \ref{poly}).

\item
For any $R$-module $M,$
$ \bigcup\limits_{a\in R} a\Gamma_{a}(M) = E_{M}(0)$ where  $E_{M}(0)$ is the envelope of the zero submodule of $M,$  (Proposition \ref{env}) and $ \bigcup\limits_{a\in R} a\Gamma_{a}(R)= \mathcal{N}(R),$ (Corollary \ref{nilrad1}).
\item For any $R$-module $M,$ $\sqrt{(0:M)}M=  \sum \limits_{a\in \sqrt{(0:M)} }a\Gamma_{a}(M),$ (Theorem \ref{sum}).
\item If $M$ is a finitely generated multiplication $R$-module, then $\sum \limits_{a\in \sqrt{(0:M)} }a\Gamma_{a}(M)=\beta(M)$ where $\beta(M)$ is the prime radical of $M,$ (Corollary \ref{mult}).

\item If $M$ is a reduced $R$-module and $\mathfrak{a}$ is the ideal of $R$ generated by $a\in R$,
 then for any $a\in R$, the $i$-th local cohomology of $M$ with respect to $\mathfrak{a}$  is given by
 $H_{\mathfrak{a}}^i(M)\cong\text{Ext}_{R}^i(R/\mathfrak{a}, M)$, (Theorem \ref{tz}).
\end{itemize}

\paragraph\noindent
 For more information about radical theory of rings, see \cite{Gard}; while for torsion theory, see \cite{Bican, Sten} among others.

\section{Reduced modules}

\paragraph\noindent
Lee and Zhou in \cite{Leezhou} introduced reduced modules. It is clear that an $R$-module $M$ is reduced if and only if for all $a\in R,$ $a\Gamma_{a}(M)=0.$ The functor $a\Gamma_{a}(-)$ therefore can be seen as a measure of how far a module is from being reduced.
\begin{defn}
Let $R$ be a ring, $M$ an $R$-module and $a\in R$. $M$ is {\emph{${a}$-reduced}} if for all $m\in M $,
 $$a^{2}m=0~ \text{ implies that}~am=0.$$
 \end{defn}
 \begin{defn}
 An $R$-module $M$ is {\emph{reduced}} if it is ${a}$-reduced for all $a\in R.$
 \end{defn}

\paragraph\noindent
It then follows that an $R$-module is (globally) reduced if and only if it is locally reduced.

\begin{prop}
Every free module defined over a reduced ring is reduced.
\end{prop}
\begin{prf}
It follows from \cite[Example 1.3]{Leezhou}.

\end{prf}
 \begin{cor}
 Any vector space is a reduced module.

 \end{cor}

  \begin{cor}\label{pro}
 A projective module defined over a reduced ring is reduced.
\end{cor}

 .
\paragraph\noindent
For an $R$-module $M,$ $a\in R$ and $k\in \Z^{+},$ we write $(0:_{M}{a}^{k})$ to denote the submodule of $M$ given by $\left\{m\in M~|~a^{k}m=0\right\}.$
\begin{prop}\label{40}
Let $M$ be an $R$-module,  $a\in R$ and ${\bf{\mathfrak{a}}}$ the ideal of $R$ generated
by $a$. The following statements are equivalent:

\begin{enumerate}
\item $M$ is $a$-reduced,

\item $a\Gamma_{a}(M)= 0$,
\item $(0:_{M}{a})~=~(0:_{M}{a}^{k})~\text{for all} ~k\in \Z^{+}$,
\item $\lim\limits_{\overrightarrow{k}}\text {Hom}_{R}(R/{\bf{\mathfrak{a}}}^{k},M)\cong \text{ Hom}_{R}(R/{\bf{\mathfrak{a}}},M)$,
\item $\Gamma_{a}(M)\cong \text{ Hom}_{R}(R/{\bf{\mathfrak{a}}},M)$,

\item $0\to \Gamma_{a}(M) \to M \to aM \to 0 $ is a short exact sequence.

\end{enumerate}
\end{prop}
\begin{prf}
\begin{itemize}

\item[$1\Ra 2$] Let $n \in a\Gamma_{a}(M).$ Then $n=am$ for some $m \in \Gamma_{a}(M).$ So, there exists $k\in \Z^{+} $ such that $a^{k}m=0.$  From $1,$ we have $am=0.$ Thus, $n=0$ and $a\Gamma_{a}(M)=0.$
\item[$2\Ra 3$] In general, $(0:_{M}{a})\subseteq (0:_{M}{a}^{k}).$ Now, let $m\in (0:_{M}{a}^{k}),$ i.e., $a^{k}m = 0.$ It follows that $m \in \Gamma_{a}(M).$ So, by $2,$ $am \in a\Gamma_{a}(M)=0 $ which implies that $m \in (0:_{M}{a}). $
\item[$3\Ra 4$] It is known that $\text{Hom}_{R}(R/{\bf{\mathfrak{a}}}^{k},M)\cong(0:_{M}{a}^{k})$ and $\lim\limits_{\overrightarrow{k}}\text {Hom}_{R}(R/{\bf{\mathfrak{a}}}^{k},M)$\\$\cong \bigcup\limits_{k\in \Z^{+}} (0:_{M}{a}^{k}),$ see \cite[page $6$]{Brodmann}.
Since by $3,$ $(0:_{M}{a}^{k})=(0:_{M}{a})$ for all $k\in \Z^{+},$ we \\have $ \lim\limits_{\overrightarrow{k}}\text {Hom}_{R}(R/{\bf{\mathfrak{a}}}^{k},M) \cong \text{ Hom}_{R}(R/{\bf{\mathfrak{a}}},M).$
\item[$4\Ra 5$] Since $\Gamma_{a}(M)\cong \lim\limits_{\overrightarrow{k}}\text {Hom}_{R}(R/{\bf{\mathfrak{a}}}^{k},M),$ $4$ implies that $\Gamma_{a}(M)\cong \text{ Hom}_{R}(R/{\bf{\mathfrak{a}}},M).$
\item[$5\Ra 6$] The $R$-module epimorphism $M \to aM$ defined by $m\mapsto am $ has kernel $(0:_{M}{a}).$ \\So, $0\to (0:_{M}{a}) \to M \to aM \to 0 $ is a short exact sequence. From $5,$ \\$\Gamma_{a}(M)\cong \text{ Hom}_{R}(R/{\bf{\mathfrak{a}}},M)\cong (0:_{M}{a})$. So, $0\to \Gamma_{a}(M) \to M \to aM \to 0$ is a short exact sequence.
\item[$6\Ra 1$] Let $a\in R$ and $m\in M $ such that $a^{2}m=0.$ Then $m\in \Gamma_{a}(M).$ From $6,$ $\Gamma_{a}(M)$ is the kernel of the epimorphism $M\to aM$ given by $m\mapsto am .$ It follows that $am=0$ which establishes $1.$
\end{itemize}
\end{prf}

\begin{prop}\label{4}
Let $M$ be an $R$-module and $\mathfrak{a}$ an ideal of $R$ generated by $a$. The following statements are equivalent:

\begin{enumerate}
\item $M$ is reduced;
\item$a\Gamma_{a}(M)= 0~\text{for all}~ a\in R;$
\item $(0:_{M}{a})=(0:_{M}{a}^{k})~\text{for all} ~a\in R, k\in \Z^{+};$
\item $\lim\limits_{\overrightarrow{k}}\text {Hom}_{R}(R/{\bf{\mathfrak{a}}}^{k},M)\cong\text{ Hom}_{R}(R/{\bf{\mathfrak{a}}},M)~\text{for all}~ a\in R;$
\item $\Gamma_{a}(M)\cong\text{ Hom}_{R}(R/{\bf{\mathfrak{a}}},M)~\text {for all } a\in R;$
\item $0\to \Gamma_{a}(M) \to M \to aM \to 0 $ is a short exact sequence for all $a\in R$.
\end{enumerate}
\end{prop}

\begin{prf}
This follows from Proposition \ref{40} and the fact that an $R$-module is reduced if and only if it is $a$-reduced for all $a \in R.$
\end{prf}

\paragraph\noindent

 Proposition \ref{40} (resp. Proposition \ref{4}) shows that being $a$-reduced (resp. reduced)
 is a categorical property, i.e., it can be expressed entirely in terms of objects and
 morphisms. Therefore, if $R$ and $S$ are rings and $F: R\text{-Mod} \to S\text{-Mod}$ is a category equivalence, then an $R$-module $M$ is $a$-reduced (resp. reduced) if and only if so is the $S$-module $F(M)$.

We now give a new characterisation of reduced rings.
\begin{cor}
Let $R$ be a ring and $\mathfrak{a}$ an ideal of $R$ generated by $a$. The following statements are equivalent:
\begin{enumerate}
\item $R$ is reduced;
\item$a\Gamma_{a}(R)= 0~\text{for all}~ a\in R;$
\item $(0:_{R}{a})=(0:_{R}{a}^{k})~\text{for all} ~a\in R, k\in \Z^{+};$
\item $\lim\limits_{\overrightarrow{k}}\text {Hom}_{R}(R/{\bf{\mathfrak{a}}}^{k},R)\cong\text{ Hom}_{R}(R/{\bf{\mathfrak{a}}},R)~\text{for all}~ a\in R;$
\item $\Gamma_{a}(R)\cong\text{ Hom}_{R}(R/{\bf{\mathfrak{a}}},R)~\text {for all } a\in R;$

\item $0\to \Gamma_{a}(R) \to R \to aR \to 0 $ is a short exact sequence for all $a\in R$.
\end{enumerate}

\end{cor}
\begin{prf}
This follows from Proposition \ref{4} and the fact that a ring $R$ is reduced if and only if the $R$-module $R$ is reduced.
\end{prf}

\paragraph\noindent

 For a given ring $R$ and $a\in R$, the submodule $a\Gamma_{a}(M)$
 is a measure of how far the $R$-module $M$ is from being $a$-reduced.

 \begin{exam}\label{p} Let $p$ be a prime number, $k\in \Z^{+}$ and $\Z_{p^{k}}$ a group of integers modulo $p^{k}.$ $\Z_{p^{k}}$ is a $ \Z$-module and $\Gamma_{p}(\Z_{p^{k}})= \Z_{p^{k}}.$  It follows that for $k=1,~\Gamma_{p}(\Z_{p})= \Z_{p}$ and hence $p\Gamma_{p}(\Z_{p})= 0,$ i.e., every $ \Z$-module  $\Z_{p}$ is $p$-reduced.
 \end{exam}

 \begin{prop}\label{propl} If $\{M_{i}\}_{i\in I}$ is a family of $R$-modules and
 $M =\prod\limits_{i\in I} M_{i}$, then $M$ is
 a reduced $(\text{resp.}~ a\text{-reduced})$ $R$-module if and only if each $M_i$ is a reduced $(\text{resp.}~a {\text{-reduced}})$ $R$-module.
\end{prop}
\begin{prf}

Follows from \cite[Example 1.3]{Leezhou}.

   \end{prf}

\begin{prop}\label{ared}
For any $R$-module $M$ and $a\in R,$ the $R$-module $M/a\Gamma_{a}(M)$ is\\ $a$-reduced.

\end{prop}
\begin{prf}
Suppose that $M/a\Gamma_{a}(M)$ is not $a$-reduced, i.e., there exists $m\in M $ and $k\in \Z^{+}$ such that $a^{k}m\in a\Gamma_{\bf{\mathfrak{a}}}(M) $ but $am\notin a\Gamma_{a}(M).$  $a^{k}m\in a\Gamma_{a}(M) $ implies that $a^{k-1}m\in \Gamma_{\bf{\mathfrak{a}}}(M). $ So, $a^{s}(a^{k-1}m)=0 $ for some $s\in \Z^{+}$ and $a^{s+k-1}m=0 .$ However, $am\notin a\Gamma_{a}(M)$ implies that $m\notin \Gamma_{a}(M)$ and as such $a^{l}m \neq 0$ for all $l\in \Z^{+},$ which is a contradiction since $a^{s+k-1}m=0. $

\end{prf}
\begin{cor}\label{rad}
For any $R$-module $M$ and $a\in R,$ $$a\Gamma_{a}(M/a\Gamma_{a}(M))=0.$$
\end{cor}
\begin{prf}
By Proposition \ref{ared}, the $R$-module $M/a\Gamma_{a}(M)$ is $a$-reduced. The desired result follows from Proposition \ref{40}.
\end{prf}

\section{Properties of the locally nilradical}\label{$*$}
\begin{paragraph}\noindent
A functor $\gamma:R\text{-Mod} \to R\text{-Mod} $ is a {\emph{preradical}} if for every $R$-homomorphism  \\$  f:M \to N,~ f(\gamma(M))\subseteq \gamma(N)$.
$\gamma $ is a {\emph{radical}} if it is a preradical and for all $M \in R$-Mod,  $\gamma(M/\gamma(M))=0.$
A radical $\gamma$ is {\emph{hereditary}} or {\emph{left exact}} if for every submodule $N$ of a module $M \in R$-Mod, $\gamma(N)= N \cap \gamma(M).$ Equivalently, if for any exact sequence $0 \to N \to M \to K $
of $R$-modules, the sequence $0 \to \gamma(N) \to \gamma(M) \to \gamma(K)$ is also exact.
\end{paragraph}

        \begin{prop}\label{nilrad}
	For any ring $R$ and $a\in R$, the  functor $$a\Gamma_{a}:R{\normalfont\text{-Mod}} \to R{\normalfont\text{-Mod}}$$
 $$M\mapsto a\Gamma_{a}(M) $$ is a   radical.
	\end{prop}
\begin{prf}
Let  $f: M \to N$ be an $R$-module homomorphism. Let $x\in f(a\Gamma_{a}(M)).$ Then $x=af(m)$ for some $m\in \Gamma_{a}(M).$ This implies that $a^{k}m=0$ for some $k \in \Z^{+}.$ So, $a^{k}f(m)=f(a^{k}m)=f(0)=0.$ This shows that $f(m)\in \Gamma_{a}(N)$ and $x= af(m)\in a\Gamma_{a}(N).$
 Hence,  $ f(a\Gamma_{a}(M))\subseteq a\Gamma_{a}(N).$  This shows that the functor $a\Gamma_{a}(-)$ is a preradical. Corollary \ref{rad}
  shows that $a\Gamma_{a}(-)$ is a radical.
\end{prf}

\paragraph\noindent
The radical $a\Gamma_{a}(-)$ is in general not left exact.
 Consider $M:=\Z_{8}$ and $N:=2\Z_{8}.$ By Example \ref{p}, if $a= 2\in \Z,$ then
 $2\Gamma_{2}(M)= 2\Z_{8}$ and $2\Gamma_{2}(N)= 4\Z_{8} \subsetneq 2\Z_{8}= N\cap 2\Gamma_{2}(M).$
However, on the subcategory of reduced $R$-modules, $a\Gamma_{a}(-)$ is a left exact radical.

\paragraph\noindent
A submodule $N$ of an $R$-module $M$ is {\it characteristic} if for all automorphisms $f$ of $M$, $f(N)\subseteq N$.

\begin{prop}\label{proj}
	Let $R$ be a ring, $a\in R$ and $M$ an $R$-module. The following statements hold.
	\begin{enumerate}
		\item $a\Gamma_{a}(R)$ is an ideal of $R.$
		\item For each $M \in R$-Mod, $a\Gamma_{a}(M)$ is a characteristic submodule of $M$ and $$a\Gamma_{a}(R)M \subseteq a\Gamma_{a}(M).$$
		\item If $M$ is  projective, then $a\Gamma_{a}(M)= a\Gamma_{a}(R)M.$
	\end{enumerate}

\end{prop}
\begin{prf}
Since $a\Gamma_{a}(-)$ is a (pre)radical, the proof follows from \cite[Proposition 1.1.3]{Bican}.

\end{prf}

\paragraph\noindent
From Proposition \ref{proj}, we can recover Corollary \ref{pro}, i.e., a projective module $M$ over a reduced ring $R$ is reduced.
 For if $R$ is reduced, then so is the module ${}_RR.$ As such, $a\Gamma_{a}(R)=0$ for all $a\in R.$ By Proposition \ref{proj}, $a\Gamma_{a}(M)=0$ for all $a\in R$ and therefore by Proposition \ref{4}, $M$ is reduced.

 \begin{prop}\label{factor}
	Let $M$ be an $R$-module, $a\in R$ and $N$ a submodule of $M.$
	\begin{enumerate}
		\item $a\Gamma_{a}(N) \subseteq N\cap a\Gamma_{a}(M)$ and $(a\Gamma_{a}(M) + N)/N \subseteq a\Gamma_{a}(M/N).$
		\item If $a\Gamma_{a}(N)= N,$ then $N \subseteq a\Gamma_{a}(M).$
		\item If $a\Gamma_{a}(M/N)=0,$ then $a\Gamma_{a}(M) \subseteq N.$
	\end{enumerate}
\end{prop}
\begin{prf}
	It follows from \cite[Proposition 1.1.1]{Bican} since $a\Gamma_{a}(-)$ is a (pre)radical.
\end{prf}
\begin{prop}\label{product}
	Let $\{M_{i}\}_{i\in I}$ be a family of $R$-modules. Then
	$$ a\Gamma_{a}\left(\bigoplus\limits_{i\in I} M_{i}\right)=\bigoplus\limits_{i\in I}a\Gamma_{a}(M_{i}) $$ and $$ a\Gamma_{a}\left(\prod\limits_{i\in I} M_{i}\right) \subseteq \prod\limits_{i\in I}a\Gamma_{a}(M_{i}).$$
\end{prop}
\begin{prf}
	It follows from \cite[Proposition 1.1.2]{Bican}.
\end{prf}

\paragraph\noindent

The radical $a\Gamma_{a}(M)$ is not idempotent. Take for instance $M:= \Z_{4}$ as a $\Z$-module.
$2\Gamma_{2}(\Z_{4})= 2\Z_{4}$  but $2\Gamma_{2}(2\Gamma_{2}(\Z_{4})) = 2\Gamma_{2}(2\Z_{4})=0.$ So $2\Gamma_{2}(2\Gamma_{2}(\Z_{4})) \neq 2\Gamma_{2}(\Z_{4}).$

\paragraph\noindent
Let $R$ be a Noetherian ring of prime characteristic $p$ and $f: R \to R$ the Frobenius ring homomorphism, i.e., $f(r)= r^{p},~ \text{for} ~r\in R.$ Let $R^{f} $ be the ring with the $R\text{-}R$ bimodule structure given by $r.s := rs$ and $s.r := sf(r)$ for $r\in R$ and $s\in R^{f}. $ $F_{R}(-):= R^{f} \otimes_{R}-$  is a right exact functor on the category  $R$-Mod and is called the {\emph{Frobenius functor}} on $R$; see \cite{Frobenius}.
\begin{thm}\label{Frob}
	Let $R$ be a Noetherian local ring of characteristic $p.$ The following statements are equivalent:\begin{enumerate}
		\item $F_{R}(-)$ is exact on $R$-Mod,
		\item $a\Gamma_{a}(-)$ is a zero functor on $R$-Mod for all $a\in R$,
		\item $R$ is a regular ring,
\item every $R$-module is reduced.
	\end{enumerate}
\end{thm}
\begin{prf}
	By \cite{Ernst}, the functor $F_{R}(-)$ is exact on $R$-Mod if and only if $R$ is a regular ring. However, by \cite[Theorem 2.16]{Rege}, $R$ is a regular ring if and only if every $R$-module is reduced, i.e., if and only if $a\Gamma_{a}(-)$ is the zero functor on $R$-Mod for all $a\in R.$
\end{prf}

\paragraph\noindent

Theorem \ref{Frob} gives a subcategory of $R$-Mod on which the Frobenius functor is exact, i.e., the subcategory of all reduced $R$-modules when $R$ is a  Noetherian local ring of characteristic $p.$ This highlights the importance of the subcategory of reduced modules over a Noetherian local ring of characteristic $p.$ They are doing to the Frobenius functor what a projective module (resp. injective module and flat module) $M$ does to the functor Hom$_{R}(M,-),$ (resp. Hom$_{R}(-,M)$ and $-\bigotimes_{R}M)$, i.e., transforming them into exact functors.

\paragraph\noindent
Since for a commutative ring $R,$ $\mathcal{N}(R)$ is the prime radical of $R$, $\mathcal{N}(R)[x]= \mathcal{N}(R[x]),$ see \cite[Theorem 10.19]{Lam}. Theorem \ref{poly} gives the local behavior of this.
\begin{thm}\label{poly}
	For any ring $R$ and $a\in R,$ $$a\Gamma_{a}(R)[x]= a\Gamma_{a}(R[x]).$$
\end{thm}
\begin{prf}

\begin{itemize}

\item[] 	Let
 $f(x)\in a\Gamma_{a}(R)[x].$ Then $f(x)= \sum \limits_{i=0}^{n}r_{i}x^{i}$ with $r_{i}\in a\Gamma_{a}(R). $ This implies that for each $r_{i},~ i\in \{0,1,2,...,n\},$ there exists $s_{i}\in \Gamma_{a}(R) $ and $k_{i} \in \Z^{+}$ such that $r_{i}=as_{i}$ and $a^{k_{i}}s_{i}=0.$ $f(x)= \sum \limits_{i=0}^{n}(as_{i})x^{i} =a\sum \limits_{i=0}^{n}s_{i}x^{i}. $ To show that $f(x)\in a\Gamma_{a}(R[x]),$ it is enough to show that $g(x)=\sum \limits_{i=0}^{n}s_{i}x^{i}\in \Gamma_{a}(R[x]) $ since $f(x)= ag(x).$
	Let $k:= \text{max}\{k_{i}\}_{i=0}^{n}.$ Then $a^{k}g(x)=\sum \limits_{i=0}^{n}a^{k}s_{i}x^{i}=0. $ Hence $g(x)\in \Gamma_{a}(R[x]) $ as required. This proves that $a\Gamma_{a}(R)[x] \subseteq a\Gamma_{a}(R[x]).$

\item[] Now, suppose that $f(x)\in a\Gamma_{a}(R[x]).$ Then $f(x)= ag(x)$ and there exists $k\in \Z^{+}$ such that $a^{k}g(x)=0.$ If $g(x)= r_{0}+...+r_{n}x^{n},$ then $f(x)= ar_{0}+...+ar_{n}x^{n}.$ We show that each $ar_{i}\in a\Gamma_{a}(R) $ for $i \in \{0,1,...,n\}.$ If $ar_{i}=0,~ ar_{i}\in a\Gamma_{a}(R). $ Suppose that $ar_{i}\neq 0.$ Then $a^{k}r_{i}=0$ for all $i\in \{0,1,...,n\}$ since $a^{k}g(x)=0.$ Then, each $r_{i}\in \Gamma_{a}(R)  $ and $ar_{i}\in a\Gamma_{a}(R)$ for all $i\in \{0, 1, \cdots, n\}$. So, $f(x)\in a\Gamma_{a}(R)[x] $  and $a\Gamma_{a}(R[x])\subseteq a\Gamma_{a}(R)[x].$	
\end{itemize}
\end{prf}

\paragraph\noindent
	$\mathcal{N}(R)=0,$ (i.e., $R$ is reduced) if and only if $a\Gamma_{a}(R)=0$ (i.e., $R$ is $a$-reduced) for each $a\in R.$   $\mathcal{N}(R)[x]= \mathcal{N}(R[x])$ and $a\Gamma_{a}(R)[x]= a\Gamma_{a}(R[x])$ for all $a\in R$.  However, $\mathcal{N}(-)$ is hereditary, see \cite[Example 3.2.12]{Gard} but the ``local" radical $a\Gamma_{a}(-)$ is not hereditary.
\paragraph\noindent
A proper submodule $N$ of an $R$-module $M$ is \emph{prime} if for all $a\in R ~\text{and}~ m\in M$, $am\in N$ implies that either $m\in N$ or $aM\subseteq N.$ A module is prime if its zero submodule is prime. A prime module is reduced. Let $\beta(M)$ denote the intersection of all prime submodules of $M.$ We call $\beta(M)$ the prime radical of $M$. Since $a\Gamma_{a}(R)\subseteq \mathcal{N}(R)$ for any ring $R$ and $a\in R$;
 and $a\Gamma_{a}(M) \subseteq \beta(M)$ for any $R$-module $M$ and $a\in R$, the locally nilradical can also be seen as a generalisation of the Baer's lower nilradical for modules.

\paragraph\noindent

Let $R$ be a ring and $a \in R.$  By \cite[Proposition 1.1.4]{Bican},
  $\mathcal{T}_{a} := \left\{M\in R\text{-Mod}  ~| ~ a\Gamma_{a}(M)= M \right \}$  is a torsion class and $ \mathcal{F}_{a} := \left \{M\in R\text{-Mod} ~ | ~ a\Gamma_{a}(M)= 0 \right \} $ is a pretorsion-free class. In general, $a$-reduced modules are not closed under extension.
$\Z_{4}$ is not a $2$-reduced $\Z$-module. However, its submodule $2\Z_{4}$ and its quotient $\Z_{4}/2\Z_{4}$ are $2$-reduced.
This shows that in general, $a$-reduced modules form a pretorsion-free class but not a torsion-free class of a torsion theory.

\section{Stratifications}

\paragraph\noindent

Let $N$ be a submodule of an $R$-module $M$. The envelope $E_{M}(N)$ of $N$ is the set
$$E_{M}(N):= \left \{am~|~ a^{k}m\in N, ~a\in R,~m\in M,~ \text{for some }k\in \mathbb{Z}^{+} \right\}.$$
 The set $E_{M}(N)$ was used by McCasland and Moore in \cite{McCasland}, Smith and Jenkins in \cite{Jenkins} and Azizi in \cite{Azizi,Azizi1} among others while studying modules and rings that satisfy the radical formula. $E_{M}(0)$ was considered as the module analogue of $\mathcal{N}(R),$ the collection of all nilpotent elements of the ring $R.$ For  any ring $R$, $E_{R}(0)= \mathcal{N}(R).$

\begin{prop}[\bf Stratification of the envelope]\label{env}
For any $R$-module M,
$$ E_{M}(0)= \bigcup\limits_{a\in R} a\Gamma_{a}(M).$$
\end{prop}

\begin{prf}
 If $m\in \bigcup\limits_{a\in R} a\Gamma_{a}(M),$ then $m = an$ for some $n\in \Gamma_{a}(M).$ This implies that $a^{k}n = 0$ for some $k\in\mathbb{Z}^{+}.$ By definition of $E_{M}(0),~m\in E_{M}(0).$ Conversely, if $m\in E_{M}(0),$ then $m = an  $ with $a^{k}n = 0$ for some $n\in M, a\in R$ and $k\in\mathbb{Z}^{+}.$ So, $n\in \Gamma_{a}(M)$ which implies that $m\in a\Gamma_{a}(M).$ Thus, $m\in \bigcup\limits_{a\in R} a\Gamma_{a}(M).$
\end{prf}

\begin{cor}[\bf Stratification of the nilradical]\label{nilrad1}
For any  ring $R,$ if $\mathcal{N}(R)$ is the collection of all nilpotent elements of $R,$ then $$ \mathcal{N}(R)= \bigcup\limits_{a\in R} a\Gamma_{a}(R)= E_{R}(0). $$
\end{cor}

\paragraph\noindent
Recall that, for an $R$-module $M$, where $R$ is a reduced ring,  the {\it torsion submodule} $t(M)$ of $M$ is the submodule
$$t(M) := \left\{ m\in M~|~am=0~\text{for some } 0 \neq a\in R\right\}.$$

\begin{prop}[\bf Stratification of the torsion submodule]\label{tors}

For any reduced module $M$ defined over a reduced ring $R$, $$t(M)= \bigcup\limits_{0 \neq a\in R} \Gamma_{a}(M).$$
\end{prop}
\begin{prf}
 If $M$ is a reduced $R$-module and $0 \neq a\in R$, then $\Gamma_{a}(M)= \left \{ m\in M~|~ am=0 \right\}.$
 So $ \Gamma_{a}(M)\subseteq t(M)$.   The reverse inclusion follows from the definitions
 of both $ \Gamma_{a}(M)$ as well as $t(M)$.
\end{prf}

\section{Comparison with other  radicals}\label{section5}

\paragraph\noindent

 A proper submodule $N$ of an $R$-module $M$ is {\emph{$a$-semiprime}} (resp. {\emph{semiprime}}) if
the $R$-module $M/N$ is $a$-reduced (resp. reduced).
 We denote by Rad$(M)$  (resp.  $S(M)$, $S_{a}(M)$)  the Jacobson radical (resp. semiprime radical,  $a$-semiprime radical)  of $M$, i.e., the intersection of all  maximal (resp. semiprime, $a$-semiprime) submodules of $M$.

\begin{prop}\label{compar}
For any $R$-module $M$ and $a\in R,$ we have the following inclusions of radical submodules of $M$: $$\bigcap\limits_{a\in R} a\Gamma_{a}(M) \subseteq a\Gamma_{a}(M)\subseteq  S_{a}(M) \subseteq S(M)\subseteq \beta(M) \subseteq {\text{Rad}}(M).$$
\end{prop}

\begin{prf}
$\bigcap\limits_{a\in R} a\Gamma_{a}(M) \subseteq a\Gamma_{a}(M)$ is trivial. If a submodule $N$
of an $R$-module $M$ is $a$-semiprime, then by definition, the module $M/N$ is $a$-reduced. By Proposition \ref{40}, $a\Gamma_{a}(M/N)= 0.$ From Proposition \ref{factor}$(3),$ we get $a\Gamma_{a}(M)\subseteq N,$ i.e., every $a$-semiprime submodule of $M$ contains the submodule $a\Gamma_{a}(M)$ of $M.$ It follows that the intersection of all $a$-semiprime submodules of $M$ contains $a\Gamma_{a}(M),$ i.e., $a\Gamma_{a}(M)\subseteq S_{a}(M).$ Since a semiprime submodule of an $R$-module $M$ is $a$-semiprime (i.e., $M/N$ reduced implies $M/N$ $a$-reduced), we have
$S_{a}(M)\subseteq S(M).$  $S(M)\subseteq \beta(M)$ is due to the fact that prime submodules are semiprime and $\beta(M) \subseteq {\text{Rad}}(M)$ follows from the fact that maximal submodules are prime.
\end{prf}

 \begin{cor}\label{compar1}
For any ring $R,$
$$\bigcap\limits_{a\in R} a\Gamma_{a}(R) \subseteq a\Gamma_{a}(R)\subseteq  S_{a}(R) \subseteq S(R)= \mathcal{N}(R) \subseteq {\text{Rad}}(R).$$
\end{cor}

\paragraph\noindent

From Proposition \ref{compar} (resp. Corollary \ref{compar1}), we can see that the radical submodule $a\Gamma_{a}(M)$ (resp. ideal $a\Gamma_{a}(R)$) is very small in comparison with other radical submodules of $M$ (resp. ideals of $R$). In addition, $a\Gamma_{a}(R)$ is a proper ideal of $R$ since it is nil and therefore the unity of $R$ cannot belong to it. 

\begin{thm}\label{sum}
	For any $R$-module $M$, $$\sqrt{(0:M)}M=\sum \limits_{a\in \sqrt{(0:M)} }a\Gamma_{a}(M).$$
\end{thm}
\begin{prf}
	Let $m \in \sqrt{(0:M)}M. $ $m= \sum \limits_{i=1}^{n}r_{i}m_{i}$ where $r_{i} \in $ $\sqrt{(0:M)},$ $m_{i} \in M$ and $n\in \Z^{+}.$ So, $r_{i}^{k_{i}}M= 0$ for all $i\in \{1,2,. ..,n\}$ and for some $k_{i} \in \Z^{+}.$ It follows that $r_{i}^{k_{i}}m_{i}= 0$ for all $i\in \{1,2,...,n\}.$ Hence, $r_{i}m_{i}\in r_{i}\Gamma_{r_{i}}(M). $ Therefore $m= \sum \limits_{i=1}^{n}r_{i}m_{i} \in  \sum \limits_{i=1}^{n}r_{i}\Gamma_{r_{i}}(M) \subseteq \sum \limits_{a\in \sqrt{(0:M)} }a\Gamma_{a}(M) $ and hence $\sqrt{(0:M)}M \subseteq \sum \limits_{a\in \sqrt{(0:M)} }a\Gamma_{a}(M) .$ Now, for any $a\in \sqrt{(0:M)} ,$ $a\Gamma_{a}(M) \subseteq \sqrt{(0:M)}M.$ So, $\sum \limits_{a\in \sqrt{(0:M)} }a\Gamma_{a}(M) \subseteq \sqrt{(0:M)}M$ which gives the reverse inclusion.
\end{prf}
\begin{cor}
	If $M$ is an $R$-module such that $(0:M)$ is a radical ideal of $R,$ then $$\sum \limits_{a\in \sqrt{(0:M)} }a\Gamma_{a}(M)=0.$$ In particular, if $a\in \sqrt{(0:M)}$, then $M$ is $a$-reduced.
\end{cor}
\begin{prf}
	By Theorem \ref{sum}, $\sqrt{(0:M)}M= \sum \limits_{a\in \sqrt{(0:M)} }a\Gamma_{a}(M). $ Since $(0:M)$ is a radical ideal of $R,$ $\sqrt{(0:M)}= (0:M).$ It follows that $\sqrt{(0:M)}M= (0:M)M= 0$ which leads to the desired result.
\end{prf}
 \begin{cor}\label{mult}
 	If the $R$-module $M$ is a finitely generated multiplication module and $\beta(M)$ is
 	 the prime radical of $M,$ then $$\sum \limits_{a\in \sqrt{(0:M)} }a\Gamma_{a}(M)= \beta(M).$$
 \end{cor}
 \begin{prf}
 	By \cite[Theorem 4]{McCasland 2}, $\sqrt{(0:M)}M= \beta(M).$ However, by Theorem \ref{sum}, $\sqrt{(0:M)}M= \sum \limits_{a\in \sqrt{(0:M)} }a\Gamma_{a}(M).$ It follows that $\sum \limits_{a\in \sqrt{(0:M)} }a\Gamma_{a}(M)= \beta(M)$ as required.
 \end{prf}

\section{Computation of local cohomology}

\paragraph\noindent
In Theorem \ref{tz}, we show that reduced modules simplify computations of local cohomology;
the usual direct limits involved in the definition of local cohomology are dropped.
\begin{thm}\label{tz}
Let $R$ be a Noetherian ring, $M$ be an $R$-module and $\mathfrak{a}$ an ideal of $R$ generated by $a\in R$. Each of the following statements holds.
\begin{enumerate}
\item If $M$ is $a$-reduced, then the $i$-th local cohomology module $H_{\mathfrak{a}}^i(M)$ is given by
$$H_{\mathfrak{a}}^i(M)\cong\text{Ext}_{R}^i(R/\mathfrak{a}, M).$$

\item If $M$ is $a$-reduced and $R/\mathfrak{a}$ is a projective $R$-module, then for all $i\geq 1$
$$H_{\mathfrak{a}}^{i}(M)=0.$$

\item If $M$ is reduced, then for all $a\in R$,
$$H_{\mathfrak{a}}^i(M)\cong\text{Ext}_{R}^i(R/\mathfrak{a}, M).$$

\end{enumerate}

\end{thm}

 \begin{prf}
 By Proposition \ref{40}, if $M$ is an $a$-reduced $R$-module, then
 $\Gamma_{a}(M)\cong\text{ Hom}_{R}(R/{\bf{\mathfrak{a}}}, M)$. The $i$-th local cohomology of $M$ which
 is the right derived functor of $\Gamma_{a}(M)$ is the $R$-module $H_{\mathfrak{a}}^i(M)\cong\text{Ext}_{R}^i(R/\mathfrak{a}, M)$. If the $R$-module $R/\mathfrak{a}$ is  projective, then it
 follows by  general theory that the module $H_{\mathfrak{a}}^i(M)$ vanishes for all $i\geq 1$.
 3 is due to the fact that, if $M$ is reduced, then it is $a$-reduced for all $a\in R$.
 \end{prf}


 ~~\\

 {\bf{Acknowledgement}}

 \paragraph\noindent
The authors would like to thank the referee for his/her valuable comments which improved this paper.
 The first author was supported by Sida bilateral programme (2015--2020) with Makerere University; Project 316: Capacity building in Mathematics and its applications and registers gratitude to Busitema University for granting her leave to undertake PhD studies at Makerere University.
 The second author wishes to thank: 1) Prof. Kobi Kremnizer of Oxford University for introducing him to torsion functors and local cohomology; and 2)
  Africa-Oxford (AfOx) initiative for supporting his visit to Oxford.

\end{document}